 \documentclass[journal,11pt,final,onecolumn]{IEEEtran}

\usepackage[pdftex]{graphicx}
\usepackage{amsfonts}
\usepackage{setspace}
\usepackage{wrapfig}
\usepackage{color}
\usepackage{cancel}
\usepackage{url}
\usepackage{enumerate}
\usepackage{cite}
\usepackage{amssymb}
\usepackage[cmex10]{amsmath}
\usepackage{subfigure}
\pdfpageattr{/Group << /S /Transparency /I true /CS /DeviceRGB>>}

\graphicspath{ {./Images/} }

\def\figurename{Fig.}

\newcommand{\field}[1]{\mathbb{#1}}                                   
\newcommand{\prn}[1]{\left(#1\right)}                                 
\newcommand{\Prn}[1]{\left[#1\right]}                                 
\newcommand{\dfe}[2]{\displaystyle\frac{d}{d#1}\left(#2\right)}       
\newcommand{\dint}[3]{\displaystyle\int_{#1}^{#2}{#3}}                
\newcommand{\abs}[1]{\left|#1\right|}                                 
\newcommand{\eltwo}[1]{\left\|#1\right\|_2}                           
\newcommand{\norm}[1]{\left\|#1\right\|}                              
\newcommand{\sign}[1]{\ \mathrm{sign}{\left(#1\right)}}								
\newcommand{\set}[1]{\displaystyle\left\{#1\right\}}									

\DeclareMathOperator*{\argmin}{arg\,min}															

\makeatletter
\newcommand{\pushright}[1]{\ifmeasuring@#1\else\omit\hfill$\displaystyle#1$\fi\ignorespaces}
\newcommand{\pushleft}[1]{\ifmeasuring@#1\else\omit$\displaystyle#1$\hfill\fi\ignorespaces}

\makeatother

\newtheorem{theorem}{Theorem}
\newtheorem{lemma}{Lemma}

\newcommand{\G}{\Gamma}
\newcommand{\Gdag}{\Gamma_{\dagger}}

\newcommand{\Gk}{\Gamma_k}

\newcommand{\Gc}{\Gamma^c}

\newcommand{\PP}{\Phi_{\Gamma}^T\Phi_{\Gamma}}

\newcommand{\GD}{\Delta}

\newcommand{\adag}{a^{\dagger}}
\newcommand{\dotadag}{\dot{a}^{\dagger}}

\usepackage{scrtime} 
\usepackage[draft,scrtime]{prelim2e}

\begin{document}
%
\title{Discrete and Continuous-time Soft-Thresholding with Dynamic Inputs}

%
%
%

\author{Aur\`ele~Balavoine,~\IEEEmembership{Student Member,~IEEE,}
		Christopher~J.~Rozell,~\IEEEmembership{Senior Member,~IEEE,}
        and~Justin~Romberg,~\IEEEmembership{Senior Member,~IEEE,}
\thanks{The authors are with the School
of Electrical and Computer Engineering, Georgia Institute of Technology, Atlanta,
GA, 30332-0250 USA}\thanks{Email: \{aurele.balavoine,crozell,jrom\}@gatech.edu.}}

\maketitle

\begin{abstract}

There exist many well-established techniques to recover sparse signals from compressed measurements with known performance guarantees in the static case. However, only a few methods have been proposed to tackle the recovery of time-varying signals, and even fewer benefit from a theoretical analysis.
In this paper, we study the capacity of the Iterative Soft-Thresholding Algorithm (ISTA) and its continuous-time analogue the Locally Competitive Algorithm (LCA) to perform this tracking in real time. ISTA is a well-known digital solver for static sparse recovery, whose iteration is a first-order discretization of the LCA differential equation. Our analysis shows that the outputs of both algorithms can track a time-varying signal while compressed measurements are streaming, even when no convergence criterion is imposed at each time step. The $\ell_2$-distance between the target signal and the outputs of both discrete- and continuous-time solvers is shown to decay to a bound that is essentially optimal. Our analyses is supported by simulations on both synthetic and real data. 


%
\end{abstract}

\makeatletter{\renewcommand*{\@makefnmark}{}
\footnotetext{This work was partially supported by NSF grants CCF-0905346 and CCF-0830456.}\makeatother}


\section{Introduction}
\label{sec:intro}

\PARstart{S}{parse} recovery studies the problem of estimating a sparse signal from a set of linear measurements, when the measurements can potentially be highly undersampled. Mathematically, the problem consists in recovering a \emph{target signal} $\adag\in\mathbb{R}^N$ from noisy linear measurements $y=\Phi\adag+\epsilon$ in $\mathbb{R}^M$, where $\epsilon$ is a noise vector and $\Phi\in\mathbb{R}^{M\times N}$ is the measurement matrix with $M\ll N$ in our scenario of interest. If the target is $S$-sparse, meaning that only $S$ coefficients in $\adag$ are non-zero, a well-studied approach consists in solving an optimization program composed of a mean-squared error term and a sparsity-inducing term (typically measured using the $\ell_1$-norm $\norm{a}_1 = \sum_i \abs{a_i}$):
\begin{equation}
\argmin\limits_{a}{\dfrac{1}{2}\left\|y-\Phi a\right\|^2_2 + \lambda \left\|a\right\|_1}.
\label{eqcont:l1}
\end{equation}
%
%
The Iterative Soft-Thresholding Algorithm (ISTA) is an algorithm that performs a discrete gradient step followed by a soft-thresholding operation at each iteration, while the Locally Competitive Algorithm (LCA) performs similar operations in a continuous-time dynamical system. Both algorithms have been shown to converge to the minimum of \eqref{eqcont:l1} when the target $\adag$ is static (see Sections \ref{ssec:lca} and \ref{ssec:ista}). 

We extend this formulation to the problem where the target signal $\adag(t)$ continuously evolves with time and generates time-varying measurements $y(t)$. While several approaches have been proposed to track a time-varying signal from compressed measurements (cf. Section \ref{ssec:review}), they currently lack performance guarantees.  
Driven by the absence of theoretical analysis for the dynamic sparse recovery problem, we study the performance of the ISTA and the LCA when their inputs vary with time.

We set our analysis in the streaming setting (i.e., data is processed as it becomes available rather than in a batch computation), since recovering a time-varying sparse signal when the measurements are streaming at a high rate or computational resources are limited is of great interest for many applications. An example of such applications is networked control systems, where a controller receives measurements from and returns a control signal to a dynamically evolving system. For instance, in \cite{nagahara_compressive_2012}, the authors propose using an extension of ISTA to obtain a sparse control signal that can be sent efficiently through a rate-limited communication channel. As a consequence, in our analysis, the two algorithms may not reach a convergence criterion before the next measurement is received. While it has been previously observed in literature that limiting the number of iterations of certain algorithms in the streaming setting could still yield good convergence results (e.g., in \cite{dai_gaussian_2011}), no analysis has been previously provided for such iteration-limited settings. 

%

  The results of our analysis are given in Section \ref{sec:tracking} and show theoretically that both ISTA and LCA can track a time-varying target and reach a favorable error bound in the streaming setting, despite not being able to fully converge for each measurement. We show that the $\ell_2$-distance between the output and the target decays with a linear rate for ISTA, while it decays with an exponential rate for the LCA. These convergence rates match those obtained in the static setting. In addition, we show that the steady-state error for each algorithm is the static error plus a term that accounts for the dynamics in the problem (in particular, the energy in the derivative of the target and the time constant of the solver). Finally, in Section \ref{sec:sim}, we present the results of simulations on both synthetic and real data that support the mathematical findings.

%
\section{Background and Related work}

\label{sec:backgrd}
%
We start by giving a review of the ISTA discrete iteration and the continuous-time LCA network, along with a summary of the results obtained in the static case. We then look at what approaches have been taken for the dynamic recovery problem.

\subsection{The Iterative Soft-Thresholding Algorithm}
\label{ssec:ista}

ISTA is one of the earliest digital algorithms developed for sparse recovery \cite{daubechies_iterative_2004}, and although it tends to converge slowly, many state-of-the-art solvers are only slight variations of its simple update rule \cite{beck_fast_2009,bioucas-dias_new_2007,figueiredo_gradient_2007,wright_sparse_2009}.
Its update rule can be seen as a discretized generalized-gradient step for \eqref{eqcont:l1}. At the $l^{\text{th}}$ iterate, the output is denoted by \mbox{$a[l]\in\mathbb{R}^N$} and the update rule is\footnote{We indicate the iterate number $l$ in brackets, to match the notation for the continuous time index, and the $n^{\text{th}}$ entry of the vector in a subscript: $a_n[l]$.}:
\begin{equation} 
\begin{split}
	a[l+1] & = T_{\lambda}\prn{ a[l] + \eta \Phi^T \prn{ y - \Phi a[l]}}.
\end{split}
\label{eq:ista1}
\end{equation}
The activation function $T_{\lambda}(\cdot)$ is the soft-thresholding function defined by:
\begin{equation}
a_n = T_{\lambda}(u_n) = \begin{cases}
0, & \abs{u_n} \leq \lambda \\
u_n-\lambda z_n, & \abs{u_n} > \lambda 
\end{cases},
\label{eqcont:thresh}
\end{equation}
where $z_n = \sign{u_n}$.
The constant $\eta$ represents the size of the gradient step, which is usually required to be in the interval $\prn{0,2\norm{\Phi^T\Phi}^{-1}}$ to ensure convergence. 
Several papers have shown that ISTA converges to the solution $a^*$ of \eqref{eqcont:l1} as $l$ goes to infinity from any initial point $a[0]$ with linear rate \cite{bredies_linear_2008, hale_fixed-point_2008}, i.e. $\exists d,\mathcal{K}$ such that $\forall l\geq0$
$$\eltwo{a[l]-a^*} \leq \mathcal{K} d^l.$$

\subsection{The Locally Competitive Algorithm}
\label{ssec:lca}

The LCA is a continuous-time system that was proposed in~\cite{rozell_sparse_2008} to solve the $\ell_1$-minimization problem and other nonsmooth sparse recovery problems. Its simple and highly-parallel architecture 
makes it amenable to analog VLSI~\cite{shapero_low_2012,shapero_configurable_2013}.
The LCA is defined by a system of coupled, nonlinear Ordinary Differential Equations (ODEs) that controls the evolution over time of a network of nodes $u_n(t)$ for $n=1, \ldots, N$:
\begin{equation} 
\begin{split}
	\tau \dot{u}(t)&=-u(t) - (\Phi^T \Phi - I) \ a(t) + \Phi^Ty\\
	a(t) &= T_{\lambda}(u(t))\\
\end{split} .
\label{eqcont:dyn}
\end{equation}
Each node or internal state variable is associated with a single dictionary element $\Phi_n$, which corresponds to a column of the matrix $\Phi$. The internal states produce output variables $a_n(t)$ for $n=1, \ldots, N$ through the activation function $T_{\lambda}(\cdot)$ applied entry-wise. The constant $\tau$ represents the time constant of the physical solver implementing the algorithm. A long lineage of research has shown that analog networks for optimization can have significant speed and power advantages over their digital counterparts~\cite{cichocki_neural_1993}.
Though $\ell_1$-minimization is our focus, previous work has shown that the LCA can be used to solve a wide class of optimization programs by varying the form of the activation function $T_{\lambda}(\cdot)$~\cite{balavoine_convergence_2012,balavoine_convergence_2013-1,charles_common_2012}.
For a wide class of activation functions, the LCA convergence is exponentially fast from any initial point $a(0)$ in the sense that $\exists d,\mathcal{K}>0$, such that $\forall t\geq0$
$$\eltwo{a(t)-a^*}\leq\mathcal{K} e^{-(1-d)t/\tau},$$
where $a^*=T_{\lambda}(u^*)$ is the solution of \eqref{eqcont:l1}~\cite{balavoine_convergence_2012}. Finally, the LCA output remains sparse when the target vector is sparse and the matrix $\Phi$ satisfies the restricted isometry property (RIP)~\cite{balavoine_convergence_2013}. In this scenario, the constant $d$ can be approximated by the RIP constant of the matrix.

To make the link between the ISTA update rule and the LCA equation clear, we rewrite the ISTA iteration \eqref{eq:ista1} as follows:
\begin{equation*} 
\begin{split}
	u[{l+1}] & = a[l] + \eta \Phi^T \prn{ y - \Phi a[l]}\\
	a[{l+1}] & =  T_{\lambda}(u[{l+1}])
\end{split} \qquad ,\forall l\geq 0.
\end{equation*}
With this formulation, it is easy to see that ISTA is a first-order discretization of the LCA dynamics (or Euler method). Using a step size $dl$ for the discretization equal to the LCA time-constant $dl = t_{l+1}-t_{l} = \tau$, the LCA ODE \eqref{eqcont:dyn} becomes:
\begin{align*}
\tau\dfrac{u[l+1]-u[l]}{\tau} &= -u[l]+a[l]+\Phi^T(y-\Phi a[l]) \\
a[l+1] & = T_{\lambda}(u[l+1]),
\end{align*}
which can be written as
\begin{align*}
u[l+1] &= a[l]+\Phi^T(y-\Phi a[l]) \\
a[l+1] & = T_{\lambda}(u[l+1]).
\end{align*}
This matches the ISTA iteration for $\eta = 1$. Because ISTA is a discrete version of the LCA dynamics, we will perform simulations on a digital computer and simulate the LCA ODE using ISTA with the appropriate parameters. This allows us to compare existing approaches to ISTA and the LCA  in Section \ref{sec:sim} by putting all the algorithms in the same framework. 

\subsection{Related work}
\label{ssec:review}


Several approaches have been proposed to tackle the problem of tracking a high-dimensional sparse signal evolving with time from a set of undersample streaming measurements. Classical methods for tracking signals include Kalman filtering and particle filtering \cite{chen_bayesian_2003}. These methods require a reliable knowledge of the underlying dynamics of the target and do not exploit any sparsity information. Some recent papers have built on these methods by incorporating a sparsity-aware criteria, either via convex relaxation \cite{vaswani_kalman_2008,charles_dynamic_2013} or greedy methods \cite{dai_gaussian_2011}, and still require a priori knowledge of the target dynamics. 

Another class of methods rely on building a probabilistic model for the evolution of the target's support and amplitudes, and use Bayesian inference techniques to estimate the next time sample \cite{sejdinovic_bayesian_2010,shahrasbi_tc-csbp:_2011,ziniel_tracking_2010}. These methods also necessitate a priori knowledge of the target's behavior to adjust several parameters. While \cite{ziniel_tracking_2010} proposes estimating the model parameters online, it only does so in the non-causal smoothing case, which can become computationally expensive as the number of parameters is large. 

Finally, the last class of methods is based on optimization. For instance, in~\cite{angelosante_compressed_2009,angelosante_online_2010}, an optimization program is set up to account for the temporal correlation in the target, and the recovery is performed in batches. In~\cite{hall_dynamical_2013}, the best dynamical model is chosen among a family of possible dynamics or parameters. In \cite{asif_dynamic_2010}, a continuation approach is used to update the estimate of the target using the solution at the previous time-step. In~\cite{zakharov_dcd-rls_2013, slavakis_adaptive_2010, babadi_sparls:_2010, chen_sparse_2009}, the optimization is solved using low-complexity iterative schemes. Unfortunately, all of the above methods lack strong theoretical convergence and accuracy guarantees in the dynamic case.
Finally, in~\cite{slavakis_generalized_2013}, a very general projection-based approach is studied. A convergence result is given but it is not clear how the necessary assumptions apply in the time-varying setting and it does not come with an accuracy result.


The ISTA and LCA belong to the class of optimization-based schemes. The two algorithms do not rely on any model of the underlying dynamics, and a minimal number of parameters need to be adjusted that are already present in the static case. In our analysis, convergence to the minimum of the objective in \eqref{eqcont:l1} or to a stopping criterion is not required. Rather, the LCA output evolves continuously with time as the input is streaming, while the standard ISTA iteration is performed as new measurements become available. This is particularly useful in scenarios where signals are streaming at very high rates or computational resources are limited. Despite this simple setting, our analysis shows that the LCA and ISTA can both track a moving target accurately and provides an analytic expression for the evolution of the $\ell_2$-distance between the output of both algorithms and the target.

\section{Tracking a time-varying input}
\label{sec:tracking}

In this section, we present the target signal model and state the two main theorems that analyze the tracking abilities of the ISTA and LCA when recovering a time-varying input $\adag(t)$.  

\subsection{Notations}

We split the indices of the states and outputs into two sets.  For $n\in\G(t)$, called the \emph{active set}, the states satisfy \mbox{$\abs{u_n(t)}>\lambda$} and the output $a_n(t)$ is strictly non-zero. For indices $n$ in the \emph{inactive set} $\G^c(t)$, the states satisfy  $\abs{u_n(t)} \leq \lambda$ and the output $a_n(t)$ is zero. We call the corresponding nodes \emph{active} and \emph{inactive} respectively. Similarly, the support of the non-zero coefficients in $\adag(t)$ is denoted as $\G_{\dagger}(t)$
We also define the set $\GD(t)$ that contains the $q$ indices with largest amplitude in $\abs{u(t)}$. While those three sets change with time as the target and algorithms evolve, for the sake of readability and when it is clear from the context, we omit the dependence on time in the notation.

We denote by $\Phi_{\mathcal{T}}$ the matrix composed of the columns of $\Phi$ indexed by the set $\mathcal{T}$, setting all the other entries to zero. Similarly, $u_\mathcal{T}$ and $a_\mathcal{T}$ refer to the elements in the original vectors indexed by $\mathcal{T}$ setting other entries to zero.

Finally, we define the sequence $\set{t_k}_{\set{k\in\field{N}}}$ of \emph{switching times} for the LCA such that $\forall t\in[t_k,t_{k+1})$, the sets $\G(t)$ and ${\Gdag}(t)$ are constant. In other word, a switch occurs if either an entry leaves or enters the support of the output $a(t)$, in which case $\G(t)$ changes, or if an entry leaves or enters the support of $\adag(t)$, in which case $\Gdag(t)$ changes.

In the following, the quantities $C_k$ refer to distinct constants that do not depend on time, and that are indexed by $k$ in order of appearance.

\subsection{Signal model}
\label{ssec:model}

\begin{figure}
\begin{minipage}[b]{\linewidth}
  \centering
  \centerline{\includegraphics[width=0.8 \linewidth]{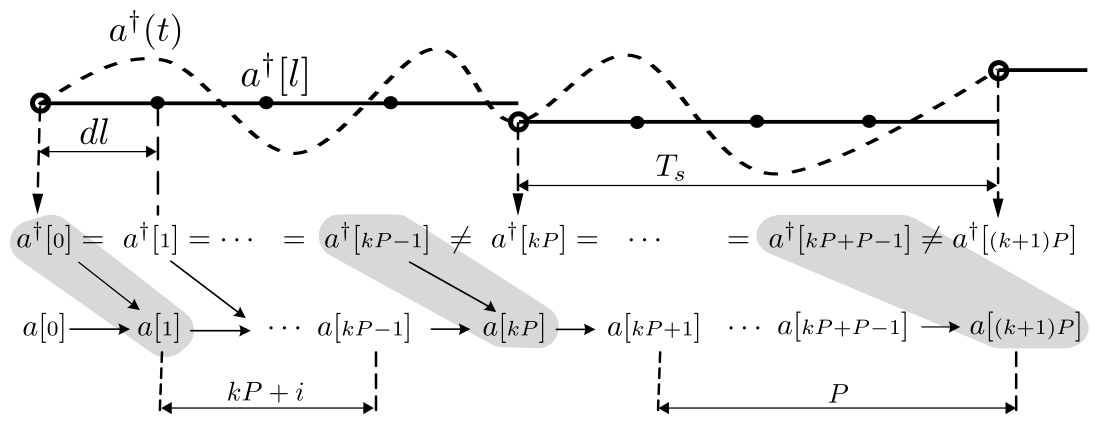}}
\end{minipage}
\caption{Measurements of the continuous-time signal $\adag(t)$ (dotted line) are continuously fed in input to the LCA. A measurement of the discrete version $\adag[l]$ (filled dots on the straight line) is received every $P^{\text{th}}$ ISTA iterates (corresponding to the empty circles). $k$ indexes the number of measurements received and $i$ indexes the number of ISTA iterations since the last measurement. The quantity of interest is $\eltwo{a[kP]-\adag[kP-1]}$ (in grey), which represents the last error before a new measurement is received.}
\label{fig:ista}
\end{figure}

In the continuous case, we assume that the underlying target signal $\adag(t)$ and the noise vector $\epsilon(t)$ evolve continuously with time. As a consequence, the input $y(t)$ to the LCA is 
\begin{equation}
 y(t) = \Phi \adag(t) + \epsilon(t), \qquad \forall t\geq0, ~t\in\mathbb{R},
\end{equation}
which is also continuous with time. 
The input to the discrete-time ISTA at iteration $l$ is a sampled version of this continuous-time signal:
\begin{equation}
y[l] = \Phi \adag[l] + \epsilon [l], \qquad \forall l\geq0, ~l\in\mathbb{Z},
\end{equation}
with sampling period $T_s$, so that $\adag[l]=\adag(l\cdot T_s)$. Letting $dl$ denote the time to perform one ISTA iteration, we state all of our results in terms of the quantities $dl$ and $P$, where $P\geq1$ is the number of ISTA iterations performed between each new measurement, so that $T_s = P\cdot dl$. Since there is no additional information about the target between two measurements, the discrete target signal $\adag[l]$ is modeled as a zero-hold version of the samples. We keep track of two indices: $k$ indexes the latest measurement received and $i$ indexes the number of ISTA iterations performed since the last measurement. Then any iteration $l$ can be expressed as $l = kP+i$ with $0\leq i \leq P-1$. This is illustrated in \figurename~\ref{fig:ista}.
With our model, iterations between new measurements that are indexed by $kP \leq l \leq (k+1)P - 1 $ satisfy
\begin{equation*}
\adag[l] = \adag[kP],
\end{equation*}
or equivalently
\begin{align}
 & \adag[kP+i] = \adag[kP], ~~\forall k\geq 0, \forall i=0,\ldots,P-1 .
\label{eq:equalenergy}
\end{align}

We focus on scenarios where the target signal is $S$-sparse (i.e., $\abs{\G_{\dagger}(t)} \leq S$, for all $t\geq 0$ and $\abs{\G_{\dagger}[l]} \leq S$, for all $l\geq 0$). 
We also require conditions on the energy in the target and its derivative. In the continuous case, our analysis requires that there exists $\mu\geq0$ such that
\begin{align}
\eltwo{\dotadag(t)} + \dfrac{1}{\tau}\eltwo{\adag(t)} & \leq \mu \qquad \forall t\geq0. \label{eqcont:targetderivenergybound}
\end{align}
 Intuitively, for a fixed $\mu$, the more energy is present in the target or the smaller the time constant $\tau$ of the solver is, the slower the target needs to vary for the LCA to track it. A consequence of condition \eqref{eqcont:targetderivenergybound} is that the energy in the time-derivative is bounded by 
 $$\eltwo{\dotadag(t)} \leq \mu, \qquad \forall t\geq0,$$
and the energy of the target itself is bounded (see Lemma~\ref{lem:targetenergybound} in Appendix \ref{app:lemmas}):
$$\eltwo{\adag(t)} \leq \beta, \qquad \forall t\geq0,$$
where $\beta = \max\set{\eltwo{\adag(0)},\tau\mu}.$ In the discrete case, these conditions simply become
\begin{equation}
\eltwo{\adag[l]} \leq \beta, \qquad \forall l\geq0,
\label{eq:bounded_energy_discrete}
\end{equation}
and, for all $k\geq0$,
\begin{align}
 \eltwo{\adag[{kP}] - \adag[kP-1]} &= \eltwo{\int_{t_{kP-1}}^{t_{kP}} \dotadag(t) dt} \nonumber \\
 &  \leq \int_{t_{kP-1}}^{t_{kP}} \eltwo{\dotadag(t)} dt \nonumber \\
 &  \leq \int_{t_{kP-1}}^{t_{kP}} \mu ~dt = \mu ~dl,
 \label{eq:targetderiv}
\end{align}
where again $dl = t_{l+1} - t_{l}$ is the time to perform one ISTA iteration. Note that because the measurement vector $y[l]$ changes every $P^{\text{th}}$ iteration, ISTA never converges to the optimum of \eqref{eqcont:l1} if $P$ is small. This approach is of great interest for scenarios where the measurements are streaming at very high rates.

Finally, we assume that the energy of the noise vector remains bounded and define $\sigma$ as:
\begin{equation}
 \eltwo{\epsilon(t)} \leq \prn{\sqrt{1+\delta}}^{-1}\sigma \qquad \forall t\geq0. \label{eqcont:noiseenergybound}
\end{equation}
Our analysis remains valid in the noise-free case, when $\epsilon = 0$. 

We assume throughout that the columns of $\Phi$ have unit norm: $\norm{\Phi_n}_2=1$ and that $\Phi$ satisfies the Restricted Isometry Property (RIP) with parameter $(S+q,\delta)$, for some $q\geq 0$. This means that for all $(S+q)$-sparse $x$ in $\mathbb{R}^N$:
\begin{equation}
(1-\delta)\eltwo{x}^2\leq\eltwo{\Phi x}^2\leq(1+\delta)\eltwo{x}^2.
\label{eq:rip}
\end{equation}
In particular $\Phi = I$ satisfies \eqref{eq:rip} with $\delta = 0$, and so direct observations of the target satisfy our requirements.

\subsection{Tracking abilities of ISTA}
\label{ssec:istares}

This section presents the tracking abilities of ISTA in the streaming setting. In the static setting, ISTA was shown in~\cite{bredies_linear_2008} to converge with a linear rate to the solution of the $\ell_1$-minimization program, for which an accuracy analysis was carried out in~\cite{zhang_sharp_2009}. Combining these two results yields an $\ell_2$-error of the following form $\forall l\geq0$:
$$\eltwo{a[l]-\adag}\leq C_0 \widetilde{c}^l + C_1,$$
where $\widetilde{c}\in(0,1)$ and $C_0,C_1\geq0$.
The constant $C_1$ represents the steady-state error $\|a^*-\adag\|$, where $a^*$ is the solution of \eqref{eqcont:l1}, and satisfies
$$C_1\leq C_2\lambda\sqrt{q}+C_3\sigma,$$
where $q\geq0$ (is typically on the order of $S$) and $C_2,C_3\geq0$. In the streaming setting, we show below that the decay of the $\ell_2$-distance between the target signal and the output of ISTA is still linear. For appropriate parameter choices, our results show that there exist $c\in(0,1)$ and two constants $C_4,C_5\geq0$ such that, for all iterations $l\geq0$,
$$\eltwo{a[l+1]-\adag[l]}\leq C_4 c^l + C_5.$$
The constant $C_5$ represents the steady-state error, after enough iterations, and satisfies
$$C_5\leq C_6\lambda\sqrt{q}+C_7\sigma+C_8\mu~dl,$$
where $C_6,C_7,C_8,\mu\geq0$. Thus, the steady-state error in the dynamic setting is simply the static error plus a term that accounts for the dynamic components of the problem.

More precisely, in the streaming setting, the ISTA $l^{\text{th}}$ iterate is
\begin{equation} 
\begin{split}
	u[{l+1}] & = a[l] + \eta \Phi^T \prn{ y[l] - \Phi a[l]}\\
	a[{l+1}] & =  T_{\lambda}(u[{l+1}])
\end{split} \qquad ,\forall l\geq 0
\label{eq:ista}
\end{equation}
and we have the following theorem.

\begin{theorem}
\label{th:istadecay}
Assume that the dictionary $\Phi$ satisfies the RIP with parameters $(S+2q,\delta)$ for some $q\geq 0$, and that the gradient step size $\eta$ in \eqref{eq:ista} satisfies
\begin{equation}
0 < \eta < \dfrac{2}{1+\delta}.
\label{eq:phitphinorm}
\end{equation}
Define \mbox{$c=\abs{\eta-1}+\delta\cdot\eta<1$}.
If the target signal satisfies conditions \eqref{eq:bounded_energy_discrete} and \eqref{eq:targetderiv}, the initial point $a[0]$ contains less than $q$ active nodes and the following two conditions hold
\begin{gather}
\eltwo{u_{\GD[0]}[0]}\leq \lambda\sqrt{q} \label{eqcont:stateenergybound2} \\
 \eta(1+\delta)\beta  + \eta\sigma \leq (1-c)\lambda\sqrt{q},\label{eq:initcond}
\end{gather}
then 
\begin{enumerate}
\item the output $a[l]$ never contains more than $q$ active nodes (i.e., $\abs{\G[l]}\leq q$ $~\forall l\geq0$); and

 \item the $\ell_2$-distance between the output and the target signal satisfies $\forall l\geq 0$, letting $i = \prn{l\mod{P}}$ ~~(i.e., \mbox{$l=kP+i$}, with $0\leq i\leq P-1$)
\begin{align}
& \qquad \qquad \qquad  \eltwo{a[l+1] - \adag[l]} \leq c^{l}\prn{\eltwo{a[1] - \adag[0]}-W} + \dfrac{c^{i+1}}{1-c^{P}}{\mu~dl} +V,
 \label{eq:istadecay}  \\
& \text{where} \nonumber \\
& \qquad \qquad \qquad V =(1-c)^{-1}\prn{ \eta\sigma + \lambda\sqrt{q}}, \label{eq:dconstant} \\
& \qquad \qquad \qquad W = \dfrac{c}{1-c^{P}}{\mu~dl} +V. \label{eq:wconstant}
\end{align}
\end{enumerate}
\end{theorem}

This theorem is proven in Appendix~\ref{app:istadecay} by induction on $l$ and gives a bound on the $\ell_2$-error between the output of ISTA and the target for every iteration $l=kP+i$. The constants $C_4$ through $C_8$ can be explicitly written as
\begin{align*}
C_4 & = \eltwo{a[1] - \adag[0]}-W \\
C_5 & \leq \underbrace{(1-c)^{-1}}_{=C_6}\lambda\sqrt{q} + \underbrace{(1-c)^{-1}\eta}_{=C_7}\sigma + \underbrace{\dfrac{c}{1-c^{P}}}_{=C_8}{\mu~dl}
\end{align*}
In addition, the result shows that at every $P^{\text{th}}$ iteration before a new measurement is received, the $\ell_2$-distance between the output and the target signal is
\begin{align*}
 \eltwo{a[kP] - \adag[kP-1]} & \leq c^{kP-1}\prn{\eltwo{a[1] - \adag[0]}-W}  + \dfrac{c^{P}}{1-c^{P}}{\mu~dl} +V.
\end{align*}
This quantity remains bounded and converges as $k$ goes to infinity towards the steady-state value
$$V+\dfrac{c^{P}}{1-c^P}\mu~dl = (1-c)^{-1}\prn{\lambda\sqrt{q}+\eta\sigma}+ \dfrac{c^{P}}{1-c^P}\mu~{dl} $$
with a linear rate of convergence. This final value is what is expected with the first term $(1-c)^{-1}\prn{\lambda\sqrt{q}+\eta\sigma}$ corresponding to the error involved with solving \eqref{eqcont:l1}. Together with the bound \eqref{eq:initcond}, these terms resemble those in Corollary 5.1 in \cite{zhang_sharp_2009} obtained for the static case. The additional term \mbox{$c^{P}(1-c^P)^{-1}\mu~{dl}$} decays with $P$ and corresponds to the error we expect from having a time-varying input. The larger the variations in the target, the larger $\mu$ is, which corresponds to a more difficult signal to track and a larger error. Conversely, the slower the target varies, the larger $P$ and as expected, the smaller the final error is. When $P$ goes to infinity, this additional term disappears.

The two conditions \eqref{eqcont:stateenergybound2} and \eqref{eq:initcond} have a similar form to the conditions of Theorem 3 in~\cite{balavoine_convergence_2013}. If there is no initial guess, then $u[0]=0$ and \eqref{eqcont:stateenergybound2} holds. Assuming $\Phi$ is a subgaussian random matrix for which the RIP constant $\delta\sim\sqrt{S/M \log\prn{N/S}}$, a similar analysis to Section IV.B of~\cite{balavoine_convergence_2013} implies that \eqref{eq:initcond} holds for a number of measurements on the order of $S\log\prn{N/S}$.

Condition \eqref{eq:phitphinorm} on the gradient step-size is better than the traditional condition that requires $0\leq\eta\leq 2 \norm{\Phi^T\Phi}^{-1}$, because the proof of Theorem~\ref{th:istadecay} uses the RIP together with the fact that the output remains $q$-sparse.

\subsection{Tracking abilities of the LCA}
\label{ssec:lcares}

In the static setting, The LCA was shown in~\cite{balavoine_convergence_2012} and \cite{balavoine_convergence_2013} to converge exponentially fast to the solution of \eqref{eqcont:l1}. For the appropriate parameter choices, the $\ell_2$-error can be expressed for all time $t\geq0$ as
$$\eltwo{a(t)-\adag}\leq C_9 e^{-v t} + C_1,$$
where $v\in(0,1)$ and $C_9\geq0$. The constant $C_1$ is again the steady-state error achieved when solving \eqref{eqcont:l1}. The theorem below shows that, in the streaming setting, the $\ell_2$-distance between the LCA output and the target signal for appropriate parameter choices can be written $\forall t\geq0$ as
$$\eltwo{a(t)-\adag(t)}\leq C_{10} e^{-v t} + C_{11}.$$
The constant $C_{11}$ represents the steady state error and satisfies
$$C_{11}\leq C_{12}\lambda\sqrt{q}+C_{13}\sigma+C_{14} \tau\mu,$$
where $C_{12},C_{13},C_{14},\mu\geq0$ with $\|\dot{a}^{\dagger}(t)\|_2\leq\mu$ for all $t\geq0$ and $\tau$ is the time constant of the LCA. As a consequence, the convergence in the dynamic setting is still exponential and the steady-state error is composed of the static error plus a term that accounts for the dynamic components in the problem.

\newpage
\begin{theorem}
\label{th:contdecay}
Assume that the dictionary $\Phi$ satisfies the RIP with parameter $(S+q,\delta)$, for some $q\geq 0$. We define the following quantity that depends on the threshold $\lambda$, the noise energy bound $\sigma$, the energy bound on the target signal $\mu$, the parameter $q$ and the RIP constant $\delta$:
\begin{equation}
 D = (1-\delta)^{-1}\prn{\tau\mu+\sigma+\lambda\sqrt{q}}. \label{eqcont:dconstant}
\end{equation}
If the target signal satisfies \eqref{eqcont:targetderivenergybound}, the initial active set $\G(0)$ contains less than $q$ active nodes, and the following two conditions hold
\begin{gather}
\eltwo{u_{\GD(0)}(0)} \leq \lambda\sqrt{q} \label{eqcont:stateenergyinit} \\
\delta \cdot \max\set{\eltwo{a(0)-\adag(0)},D} + \beta + \sigma \leq \lambda\sqrt{q}, \label{eqcont:decaycond}
\end{gather}
then 
\begin{enumerate}
\item the output $a(t)$ never contains more than $q$ active nodes (i.e., $\abs{\G(t)}\leq q$ $~\forall t\geq0$ ); and
%
%
%
 \item the $\ell_2$-distance between the LCA output and the target signal satisfies $\forall t\geq 0$
\begin{equation}
 \begin{split}
 \eltwo{a(t) - \adag(t)} &\leq e^{-(1-\delta)t/\tau}\eltwo{a(0)-\adag(0)}  +\prn{1-e^{-(1-\delta)t/\tau}} D. 
 \end{split} \label{eqcont:outputdecay}
\end{equation}
\end{enumerate}
\end{theorem}

This theorem is proven by induction on the sequence of switching times in Appendix \ref{app:contdecay} and shows that the $\ell_2$-distance between the LCA output and the target signal converges exponentially fast towards its final value $D$ in \eqref{eq:dconstant}. This bound is again what is expected for the problem. The first term $(1-\delta)^{-1}\prn{\lambda\sqrt{q}+\sigma}$ corresponds to the expected error when solving \eqref{eqcont:l1} as shown in \cite{zhang_sharp_2009}, while the additional term \mbox{$(1-\delta)^{-1}\tau\mu$} corresponds to the error associated with recovering a time-varying signal. Again, we see that the error increases with $\mu$, the energy of the variations in the target. Conversely, the error decreases with a decreasing $\tau$, corresponding to a faster solver. 

The steady state value $D$ for the LCA is equal to the steady-state value obtained for ISTA when its parameters are chosen appropriately. When ISTA is considered as the discretization of the LCA, then the discretization step $dl$ for the Euler method is equal to $\tau$, $\eta = 1$ in \eqref{eq:ista}, $P=1$ and so $c=\delta$ as discussed in Section \ref{ssec:lca}. Then the steady-state value for ISTA becomes:
$$V+\dfrac{c^P}{1-c^P}\mu~dl\leq(1-\delta)^{-1}\prn{\tau\mu + \sigma + \lambda\sqrt{q}}.$$

The initial conditions \eqref{eqcont:stateenergyinit} and \eqref{eqcont:decaycond} are similar to the initial conditions of Theorem 3 in \cite{balavoine_convergence_2013}. In particular, \eqref{eqcont:stateenergyinit} holds when there is no initial guess (i.e., $u(0)=0$) and a similar analysis to that of Section IV.B in \cite{balavoine_convergence_2013} shows that the number of measurements for \eqref{eqcont:decaycond} to hold is on the order of $S\log\prn{N/S}$ for classic CS matrices for which $\delta\sim \sqrt{S/M\log\prn{N/S}}$.

A careful look at the proof in Appendix \ref{app:contdecay} reveals that the following condition is necessary for our analysis:
\begin{align*}
& \eltwo{\dotadag(t)} + \dfrac{1}{\tau}\eltwo{\adag_{\Gc}(t)} \leq \mu \qquad \forall t\geq0.
\end{align*}
This is less restrictive than \eqref{eqcont:targetderivenergybound} since as the LCA evolves, the output gets closer to the target signal and the energy in $\|\adag_{\Gc}(t)\|_2=\|\adag_{\Gc}(t) - a_{\Gc}^{}(t)\|_2$ decreases (the equality comes from the fact that $\Gc(t)$ is the inactive set of $a(t)$ by definition). However, because the set $\Gc(t)$ changes with time and the sequence of active sets is not known in advance, this condition is difficult to verify a priori in practice.

 As a final remark, the convergence rate of both continuous and discrete algorithms depend on the RIP constant $\delta$ of the matrix $\Phi$. However, the condition on the RIP parameters $(S+2q,\delta)$ is stronger in Theorem \eqref{th:istadecay} (vs. $(S+q,\delta)$ in Theorem \eqref{th:contdecay}). This is because the ISTA is a discrete-time algorithm and the set of active elements $\G[l+1]$ may differ by as much as $q$ elements from the previous active set $\G[l]$, while the changes are continuous in the case of the LCA.

\newpage
\section{Simulations}
\label{sec:sim}
In this section, we provide simulations that illustrate the previous theoretical results\footnote{Matlab code running the experiments in this section can be downloaded at \url{http://users.ece.gatech.edu/~abalavoine3/code}}.

\subsection{Synthetic data}

We generate a sequence of sparse vectors $\set{\adag[l]}_l$ for $l=1,\ldots,40$ of length $N=512$ with sparsity $S=40$ and with up to $20$ entries that are allowed to leave or enter the support in the following manner. First, we generate the sequence of time-varying amplitudes for the nonzero entries. For $l=1$, $S=40$ amplitudes are drawn from a normal Gaussian distribution and normalized to have norm $\beta$. With these initial samples calculated, 39 consecutive time samples are generated as
$$\alpha[l+1] = \sqrt{\dfrac{\beta^2-\mu^2}{\beta^2}}\alpha[l] + \dfrac{\mu}{\sqrt{S}} v[l],$$
where $v[l]$ is a vector in $\mathbb{R}^S$ with amplitudes drawn from a standard Gaussian distribution. Each sample $\alpha[l]$ in the sequence has energy equal to $\beta$ in expectation and innovation with energy proportional to $\mu$. A set of $30$ random indices is selected and $30$ of the sequences in $\set{\alpha[l]}_{l=1,\ldots,40}$ are assigned to the corresponding indices in $\set{\adag[l]}_{l=1,\ldots,40}$. 
Finally, each of the $10$ remaining sequences in $\set{\alpha[l]}_{l=1,\ldots,40}$ is assigned to two random indices in $\set{\adag[l]}_{l=1,\ldots,40}$, and a sinusoidal envelope is generated that controls when each index is active.  The first index is set to the product of the envelope with the sequence $\alpha[l]$ when the envelope is positive and zero otherwise.  Likewise, the second index is set to the product of the envelope and sequence when the envelope is negative and zero otherwise.  This configuration results in a smooth change in support, while also enforcing that exactly $S=40$ elements are nonzero in $\adag[l]$ at all times.

\begin{figure}
\centering
	\subfigure[Effect of the parameter $P$.]{
  \includegraphics[width=0.465\linewidth]{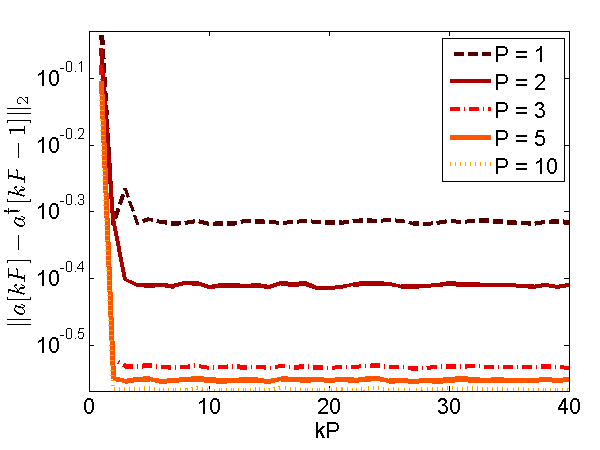}
	\label{fig:varyP} }
	\subfigure[Effect of the parameter $\mu$.]{
  \includegraphics[width=0.465\linewidth]{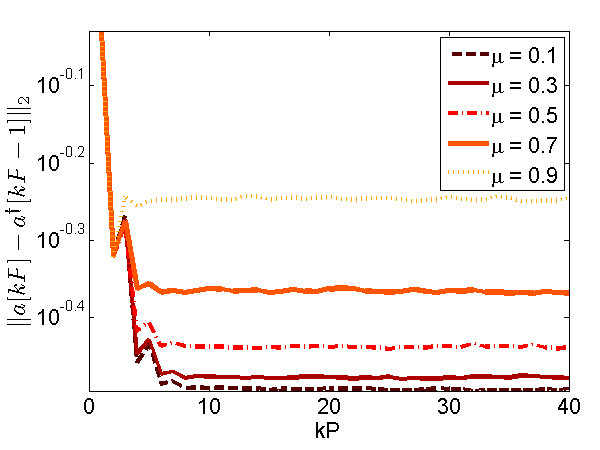}
	\label{fig:varyMu} }
	\caption{$\ell_2$-distance between the target and the output at the $P^{\text{th}}$ ISTA iterations.}
	\label{fig:vary}
\end{figure}

\begin{figure}
\begin{minipage}[b]{\linewidth}
  \centering
  \centerline{\includegraphics[width=0.465\linewidth]{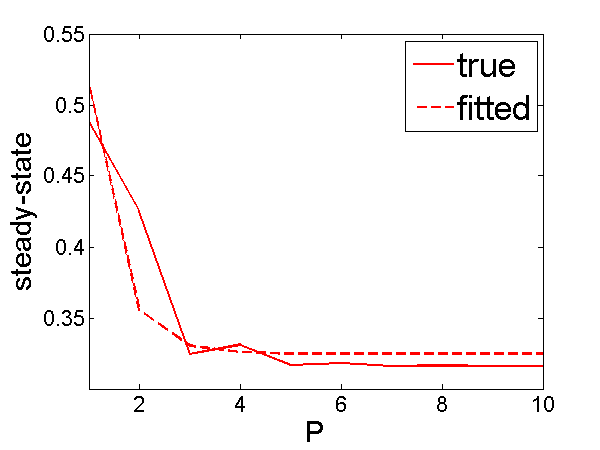}}
\end{minipage}
\caption{Steady-state of the ISTA as it recovers synthetic time-varying signals and best fit of the function $c^{P}\prn{1-c^P}^{-1}\mu~dl + V$ as a function of $P$.}
\label{fig:steady}
\end{figure}

When they are not varied, we choose $\eta = 1$, $\mu = 0.8$ and $P=1$. The measurement matrix $\Phi$ is $256\times 512$ with entries drawn from a random normal Gaussian distribution and columns normalized to 1. A Gaussian white noise with standard deviation $0.3\eltwo{\Phi \adag[0]}/\sqrt{M}$ is added to the measurements. This corresponds to a moderate level of noise. In \figurename~\ref{fig:vary}, we plot the average over 1000 such trials of the $\ell_2$-error $\eltwo{a[kP]-\adag[kP-1]}$. We observe that the curves tend to a final value that matches the steady-state behavior $c^{P}\prn{1-c^P}^{-1}\mu~dl$ predicted by theorem \ref{th:istadecay} as $k$ goes to infinity. A higher value of $P$ or a lower value of $\mu$ yields a lower steady-state. To verify that the theorem captures the correct behavior, we also plot on \figurename~\ref{fig:steady} the steady-state value as a function of $P$ for $\mu=0.8$ and find the parameters $c$ and $V$ so that the predicted function $c^{P}\prn{1-c^P}^{-1}\mu~dl + V$ best fits the data. As the figure shows, this function seems to represent the steady-state behavior well.

\begin{figure}
\begin{minipage}[b]{\linewidth}
  \centering
  \centerline{\includegraphics[width=0.465 \linewidth]{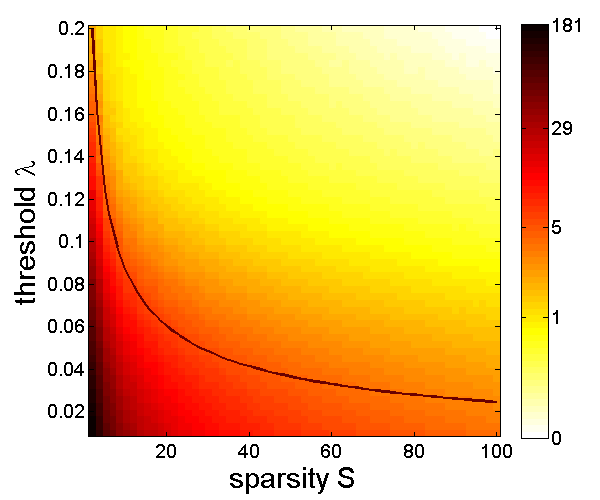}}
\end{minipage}
\caption{Ratio of the maximum number of non-zero elements $q$ over the sparsity level $S$ for several values of $\lambda$ and $S$ averaged over 100 trials.  The curve superimposed represents the best fit of the form $\lambda = C/\sqrt{S}$ of the curve for which $q=4$ is constant.}
\label{fig:qratio}
\end{figure}

Next, we vary the threshold $\lambda$ and the sparsity level $S$, and for each pair generate 10 time samples of $\adag[k]$ and associated measurements $y[k]$ in the same fashion as above. We run ISTA for $P=5$ iterations per measurement. In \figurename~\ref{fig:qratio}, we plot the average over 100 such trials of the ratio of the maximum number of non-zero elements $q$ in $a[l]$ over the sparsity level $S$. The plot shows that the maximum number of non-zero elements remains small ($q$ is mostly contained between $1S$ and $10S$), which matches the theorems' prediction. The curve superimposed shows that $\lambda$ and $S$ obey a relationship of the form $\lambda = C/\sqrt{S}$, as expected from conditions \eqref{eq:initcond} and \eqref{eqcont:decaycond}.


\subsection{Real Data}

\begin{figure*}
\begin{minipage}[b]{\linewidth}
  \centering
  \subfigure[Average rMSE]{
	  \includegraphics[width=0.465\linewidth]{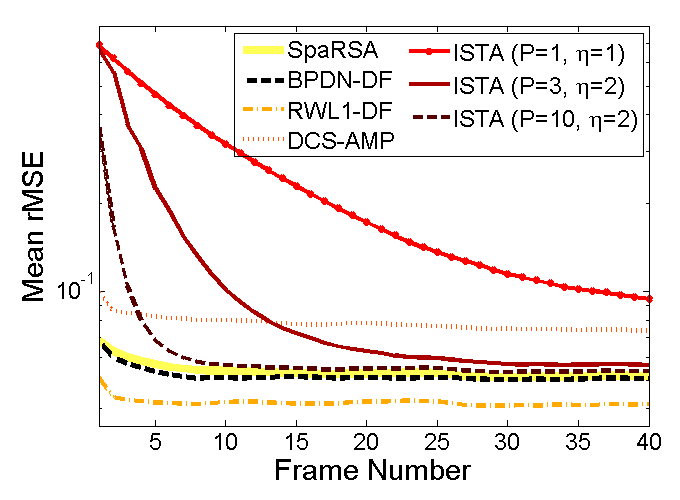}
	  \label{fig:rMSE}}
  \subfigure[Average number of products by $\Phi$ and $\Phi^T$]{
	  \includegraphics[width=0.465\linewidth]{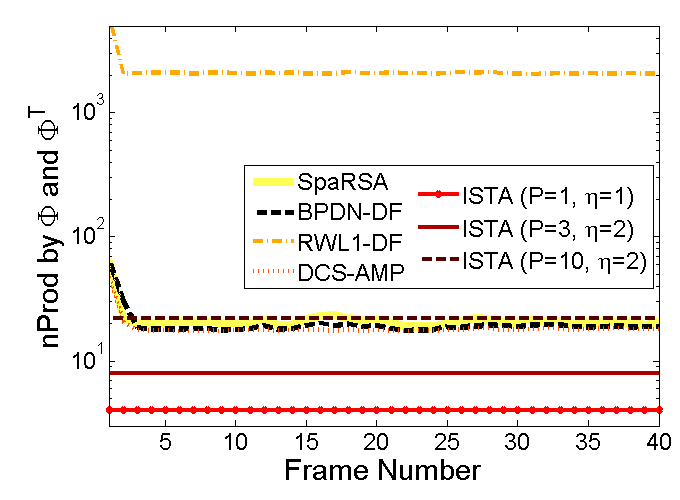}
	  \label{fig:nProd} }
	  \caption{Results of the experiment to recover the wavelet coefficients averaged over 100 random video sequences of 40 consecutive frames.}
	  \label{fig:comp}
\end{minipage}
\end{figure*}

Finally, we test how ISTA performs in the streaming setting on real data and compare it against SpaRSA, a state-of-the-art LASSO solver \cite{wright_sparse_2009}, BPDN-DF that adds a time-dependent regularization between frames, RWL1-DF which additionally performs reweighting at each iteration \cite{charles_dynamic_2013} and DCS-AMP that uses a probabilistic model to describe the target's evolution \cite{ziniel_tracking_2010}. We use 13 videos of natural scenes and randomly select 100 random sequences of 40 frames~\footnote{The videos used can be downloaded at \url{http://trace.eas.asu.edu/yuv/}}. Since natural images are sparse in wavelet basis, following the work in \cite{romberg_imaging_2008}, we use $\Phi = AB$ as the measurement matrix, where $A$ consists of $M=0.25N$ random rows of a noiselet matrix and $B$ is a dual-tree discrete wavelet transform (DT-DWT) \cite{selesnick_dual-tree_2005}. SpaRSA and BPDN-DF are given the previous solution as a warm start for the following frame. The results obtained for ISTA with $P=1$ and $\eta=1$ simulate the LCA ODE's. In \figurename~\ref{fig:rMSE}, we plot the relative mean-squared error, defined by \mbox{$rMSE[k] = \dfrac{\eltwo{a[k]-\adag[k]}}{\eltwo{\adag[k]}}$}, and the number of products by the matrix $\Phi$ or its transpose in \figurename~\ref{fig:nProd}, averaged over the 100 trials. We use the number of products by $\Phi$ and $\Phi^T$ rather than CPU time because it is a less arbitrary measure of the computational complexity for each algorithm. Figure~\ref{fig:rMSE} shows that the average rMSE values for ISTA with $P=3$ and $P=10$ reach a similar value to the rMSEs for SpaRSA and BPDN-DF after about 19 and 6 frames respectively. The average rMSE for ISTA with $P=1$ converges much slower. The rMSE for RWL1-DF is much lower, due to the additional reweighting steps, however, the complexity for this method is much larger than the other four. While the complexity of ISTA for $P=10$ is similar to the ones for SpaRSA, BPDN-DF and DCS-AMP, the complexity for $P=3$ and $P=1$ can be much smaller than any of the other approaches. In addition, ISTA only requires that two parameters be adjusted $\lambda$ and $\eta$, while DCS-AMP has approximately 10. 

%
%

\section{Conclusions}
\label{sec:ccl}

 The two main results of this paper provided an analysis for both the discrete-time ISTA and the continuous-time LCA in the context of the online recovery of a time-varying signal from streaming compressed measurements. Contrary to traditional approaches that require convergence or a stopping criterion before moving on to recover the next frame, we studied the outputs of the LCA and ISTA when their evolution rate is dictated by the rate at which measurements are streaming. Our results show that both algorithms are able to track a time-varying target in this setting. We believe that the techniques developed in this paper will prove useful for the analysis of other algorithms that currently lack theoretical analysis, in particular iterative-thresholding schemes that extend the classic ISTA. The results presented are particularly relevant for applications where the recovery needs to be performed in real-time and the signals are sampled at high rates or the computational resources are limited. While the simulations of the LCA ODE's with the appropriate parameters suggest that its convergence is slow, one would expect its time constant to be actually much faster in an actual analog system, so that the actual behavior would be closer to the one for ISTA with $P=10$ (corresponding to an analog constant 10 times smaller than the digital equivalent). If this is the case, an analog implementation of the LCA has the potential to lead to a real-time and low-power solver for such challenging situations.


\appendices


\newpage
\section{Differential Equations}
\label{app:ode}

We will use the following results on differential equations freely in the proofs.
\begin{lemma}
\label{lem:odesol}
If $x(\cdot):\field{R}^+\rightarrow\field{R}$ satisfies
\begin{align*}
 & \dfrac{d x(t)}{dt} \leqq -ax(t) + F(x(t),t), \\
 & x(0) = x_0,
\end{align*}
then $\forall t\geq 0$:
\begin{equation}
x(t) \leqq e^{-at}x_0 + e^{-at} \int_0^t e^{as} F(x(s),s) ds.
\label{eqcont:odesol}
\end{equation}
\end{lemma}

\begin{IEEEproof}
To show this, we compute $\forall t\geq 0$
\begin{align*}
 \dfrac{d}{dt}\prn{e^{at}x(t)} & = ae^{at}x(t)+e^{at}\dfe{t}{x(t)}\\
 & \leqq ae^{at}x(t)+e^{at}\prn{-ax(t) + F(x(t),t)}\\
 & \leqq e^{at}F(x(t),t).
\end{align*}
Integrating on both sides from $0$ to $t$ (using the positivity of the integral) yields
\begin{align*}
 e^{at}x(t) - x(0) & \leqq \int_0^t e^{as} F(x(s),s) ds.
\end{align*}
Thus:
\begin{align*}
 x(t) \leqq e^{-at}x_0 + e^{-at} \int_0^t e^{as} F(x(s),s) ds, \tag*{\IEEEQED}
\end{align*}
\let\IEEEQED\relax
\end{IEEEproof}


\section{Lemmas}
\label{app:lemmas}

The following lemmas are also used in the proofs. The first one is a consequence of the RIP property in \eqref{eq:rip}.

\begin{lemma}
\label{lem:ripcsq}
If $\Phi$ satisfies the RIP with parameters $(S+q,\delta)$, the set $\Gamma_1$ contains less than $q$ indices and the set $\Gamma_2$ contains less than $S$ indices, then $\forall x\in\mathbb{R}^N$ supported on $\Gamma_1 \cup \Gamma_2$ and $\forall y\in\mathbb{R}^N$, the following holds:
\begin{align*}
 (1-\delta)\eltwo{x}^2 &\leq \eltwo{\Phi x}^2 \leq (1+\delta) \eltwo{x}^2, \\
 \eltwo{\Phi_{\Gamma_1}^T\Phi_{\prn{\Gamma_1^c \cap \Gamma_2}} x} &\leq \delta\eltwo{x}, \\
 \eltwo{\prn{I_{\Gamma_1} - \Phi_{\Gamma_1}^T\Phi_{\prn{\Gamma_1\cup\Gamma_2}}}x}&\leq\delta\eltwo{x},\\
 \eltwo{\Phi_{\Gamma_1}^T y} &\leq \sqrt{1+\delta}\eltwo{y}.
\end{align*}
\end{lemma}
\begin{IEEEproof}
The RIP implies that the singular values of $\Phi_{\prn{\Gamma_1\cup\Gamma_2}}$ are contained between $\sqrt{1-\delta}$ and $\sqrt{1+\delta}$, which entails the first inequality. The singular values of $\Phi_{\Gamma_1}^T$ are contained between $0$ and $\sqrt{1+\delta}$, which implies the last inequality. Finally, the eigenvalues of $\Phi_{\prn{\Gamma_1\cup\Gamma_2}}^T\Phi_{\prn{\Gamma_1\cup\Gamma_2}}$ are contained between $(1-\delta)$ and $(1+\delta)$. For the second and third inequalities, it suffices to observe that $\Phi_{\Gamma_1}^T\Phi_{\prn{\Gamma_1^c \cap \Gamma_2}}$ and $\prn{I_{\Gamma_1} - \Phi_{\Gamma_1}^T\Phi_{\prn{\Gamma_1\cup\Gamma_2}}}$ are submatrices of $\prn{I_{\prn{\Gamma_1\cup\Gamma_2}} - \Phi_{\prn{\Gamma_1\cup\Gamma_2}}^T\Phi_{\prn{\Gamma_1\cup\Gamma_2}}}$. As a consequence, their operator norm is bounded by the operator norm of the bigger matrix, which is:
\begin{align*}
 & \norm{I_{\prn{\Gamma_1\cup\Gamma_2}} - \Phi_{\prn{\Gamma_1\cup\Gamma_2}}^T\Phi_{\prn{\Gamma_1\cup\Gamma_2}}} \leq \max\set{(1+\delta)-1,1-(1-\delta)} = \delta. 
\end{align*}
\end{IEEEproof}

The following lemma gives a bound on the energy of the target when its time-derivative satisfies \eqref{eqcont:targetderivenergybound}.
\begin{lemma}
\label{lem:targetenergybound}
If the target signal $\adag(t)$ is continuous and satisfies \eqref{eqcont:targetderivenergybound} for all $t\geq 0$
then, $\forall t\geq 0$
\begin{align*}
 \eltwo{\adag(t)} &\leq e^{-t/\tau}\prn{\eltwo{\adag(0)}-\tau\mu} + \tau\mu \\
 & \leq \underbrace{\max\set{\eltwo{\adag(0)},\tau\mu}}_{=\beta}.
 \end{align*}
\end{lemma}

\begin{IEEEproof}
It suffices to notice that 
$$\dfrac{d}{dt}\left(\eltwo{x(t)}\right) = \dfrac{\frac{d}{dt}\prn{\eltwo{x(t)}^2}}{2\eltwo{x(t)}}  = \dfrac{x(t)^T\dot{x}(t)}{\eltwo{x(t)}}\leq \eltwo{\dot{x}(t)},$$
where the last inequality comes from the Cauchy-Schwartz inequality. Thus, \eqref{eqcont:targetderivenergybound} implies
$$\dfrac{d}{dt}\left(\eltwo{\adag(t)}\right)\leq -\dfrac{1}{\tau}\eltwo{\adag(t)} + \mu.$$
Since $\adag(t)$ is continuous, we can apply Lemma \ref{lem:odesol} to obtain the first inequality. The second inequality immediately follows from the monotonicity of the exponential.
\end{IEEEproof}

The final lemma shows that a certain bound on the energy of the $q$ entries with largest absolute value $\GD$ in $u(t)$ imposes that at most $q$ nodes are active.
\begin{lemma}
\label{lem:inac}
 If $\GD$ contains the indices of the $q$ entries with largest absolute values in $u(t)$ and
 $$\eltwo{u_{\GD}(t)} \leq \lambda\sqrt{q},$$
then the active set $\G$ corresponding to the non-zero elements in $a(t) = T_{\lambda}(u(t))$ is  subset of $\GD$ and contains less than $q$ indices, i.e. $\G\subset\GD$ and $\abs{\G}\leq q$.
\end{lemma}

\begin{IEEEproof}
 Since $\GD$ contains the $q$ nodes with largest absolute values in $u(t)$, then $\forall j\in\GD^c$, we have
 $$\abs{u_j(t)} \leq \eltwo{u_{\GD}(t)}/\sqrt{q} \leq\lambda.$$
As a consequence, nodes in $\GD^c$ are below threshold, which shows that only the nodes in $\GD$ can be non-zero in $a(t)$. 
\end{IEEEproof}


\section{Proof of Theorem~\ref{th:istadecay}}
\label{app:istadecay}
\begin{IEEEproof}
We start by checking that $c<1$. By \eqref{eq:phitphinorm}, we have
\begin{align*}
 \eta(1+\delta) < 2 \text{ and } \delta<1 & \Rightarrow \eta\delta < \eta < 2-\eta\delta \\
 & \Rightarrow -1+\eta\delta < \eta - 1 < 1-\eta\delta \\
 & \Rightarrow \abs{\eta - 1} < 1-\eta\delta \\
 & \Rightarrow c<1.
\end{align*}
In a first step, we show that 
\begin{equation}
\eltwo{u[l]_{\GD[l]}}\leq\lambda\sqrt{q}, \qquad \forall l\geq0
\label{eq:statebound}
\end{equation}
by induction on $l$. If so, by Lemma \ref{lem:inac}, the active set contains less than $q$ elements and 1) of the Theorem is proven. By \eqref{eqcont:stateenergybound2}, \eqref{eq:statebound} holds for $l=0$. Next, we assume that \eqref{eq:statebound} holds for some $l\geq 0$. By Lemma \ref{lem:inac}, we have $\G[l]\subset\GD[l]$ and $\abs{\G[l]}\leq q$. As a consequence, the set 
$$J=J[l+1]:=\GD[l+1]\cup\G[l]\cup\Gdag[l]$$
contains less than $S+2q$ indices. Using the RIP of $\Phi$, the eigenvalues of the matrix $\Phi_J^T\Phi_J$ are contained between $(1-\delta)$ and $(1+\delta)$ and
\begin{align*}
 \norm{\eta\Phi_{J}^T\Phi_J-I_{J}} & = \max\set{\abs{\eta(1+\delta) - 1},\abs{\eta(1-\delta)-1}} \\
 & = \abs{\eta-1} + \eta\delta  ~~= c.
\end{align*}
In addition, the form of the activation function \eqref{eqcont:thresh} implies that $\eltwo{a[l]}\leq\eltwo{u_{\G[l]}[l]}\leq\eltwo{u_{\GD[l]}[l]}\leq\lambda\sqrt{q}$.
Using \eqref{eq:bounded_energy_discrete} and hypothesis \eqref{eq:initcond}, we get
\begin{align*}
 \eltwo{u_{J}[l+1]} & = \left\|  \eta\Phi_{J}^T\Phi\prn{\adag[l]-a[l]} + a_{J}[l] + \eta\Phi_{J}^T\epsilon[l]  \right\|_2 \\
 & \leq \eltwo{ \Prn{\eta\Phi_{J}^T\Phi_J-I_{J}} a[l] } + \eltwo{\eta\Phi_{J}^T\Phi_J \adag[l]} +\eta\eltwo{\Phi_{J}^T\epsilon[l]}  \\
 &  \leq c\eltwo{a[l]} + \eta(1+\delta)\eltwo{\adag[l]} +\eta\sqrt{1+\delta}\eltwo{\epsilon[l]} \\
 &  \leq c\lambda\sqrt{q} + \eta(1+\delta)\beta + \eta\sigma \\
 &  \leq \lambda\sqrt{q}. 
\end{align*}
Since $\GD[l+1]\subset J[l+1]$, the induction hypothesis \eqref{eq:statebound} holds at $l+1$. As a consequence, we proved 1) of the theorem and also the stronger result $\eltwo{u[l]_{J[l]}}\leq\lambda\sqrt{q}$, $~\forall l\geq1$, which will be used in the following.

Next, we show by induction on $l$ that \eqref{eq:istadecay} holds $\forall l\geq0$. It obviously holds for $l=0$. Next, assume that for some $l\geq 0$, \eqref{eq:istadecay} holds. There exist a unique $k\geq0$ and a unique $0\leq i\leq P-1$ such that $l=kP+i$. In the previous part of the proof, we showed that $\eltwo{u_{J'}[l+2]}\leq\lambda\sqrt{q}$, where $J'=J[l+2]=\GD[l+2]\cup\G[l+1]\cup\Gdag[l+1]$ and that $J'$ contains less than $S+2q$ indices. As a consequence, we can use the RIP of $\Phi_{J'}^T\Phi_{J'}$ and get
\begin{align*}
 \eltwo{a[{l+2}]-\adag[{l+1}]}  & \leq \eltwo{a[{l+2}]-u_{J'}[l+2]}  + \eltwo{u_{J'}[l+2]-\adag[{l+1}]} \\
 %
 %
 & \leq \eltwo{u_{J'}[l+2]}+ \eltwo{u_{J'}[l+2]-\adag[{l+1}]}  \\
 &  \leq \lambda\sqrt{q} + \left\| \eta\Phi_{J'}^T\epsilon[l+1] + \prn{\eta\Phi_{J'}^T\Phi - I_{J'}} \prn{\adag[l+1] - a[l+1]} \right\|_2  \\
 &  \leq \lambda\sqrt{q} + \eta\sigma + c\eltwo{\adag[l+1] - a[l+1]}\\
 &  \leq \lambda\sqrt{q} + \eta\sigma + c\eltwo{a[l+1]-\adag[l]}  + c\eltwo{\adag[l]-\adag[l+1]} .
\end{align*}
We use the induction hypothesis \eqref{eq:istadecay} at $l$ and split the analysis in two cases:

\textbf{First case:}
When $i=P-1$, $l= (k+1)P -1$ and \eqref{eq:targetderiv} yields $\eltwo{\adag[l+1]-\adag[{l}]} \leq \mu~dl$. Thus,
\begin{align*}
 \eltwo{a[{l+2}]-\adag[{l+1}]} & \leq c \bigg( c^{l}\Prn{\eltwo{a[1] - \adag[0]} -W} + \dfrac{c^P}{1-c^P}\mu{~dl} + V \bigg) + c\mu{~dl} + \lambda\sqrt{q} + \eta\sigma \\
 & \leq c^{l+1}\Prn{\eltwo{a[1] - \adag[0]} -W} + \dfrac{c^{P+1}}{1-c^P}\mu{~dl} + cV + c\mu{~dl} + \lambda\sqrt{q} + \eta\sigma \\
 %
 %
 &  \leq c^{l+1}\Prn{\eltwo{a[1] - \adag[0]} -W} + \dfrac{c}{1-c^P}\mu{~dl} + V.
\end{align*}
So the induction hypothesis \eqref{eq:istadecay} holds for $l+1=(k+1)P$.

\textbf{Second case:}
When $0\leq i\leq P-2$, \eqref{eq:equalenergy} yields $\eltwo{\adag[l+1]-\adag[{l}]} =0$ and so
\begin{align*}
 \eltwo{a[{l+2}]-\adag[{l+1}]} & \leq c \bigg( c^{l}\prn{\eltwo{a[1] - \adag[0]} -W} + \dfrac{c^{i+1}}{1-c^P}\mu{~dl} +V \bigg) + \lambda\sqrt{q} + \eta\sigma \\
 &  \leq c^{l+1}\Prn{\eltwo{a[1] - \adag[0]} -W} + \dfrac{c^{i+2}}{1-c^P}\mu{~dl} + cV + \lambda\sqrt{q} + \eta\sigma \\
 &  \leq c^{l+1}\Prn{\eltwo{a[1] - \adag[0]} -W} + \dfrac{c^{i+2}}{1-c^P}\mu{~dl} + V.
\end{align*}
Since $l+1 = kP + (i+1)$, with $1\leq i+1 \leq P-1$, this proves the induction hypothesis \eqref{eq:istadecay} in the second case and finishes the proof.
\end{IEEEproof}


\section{Proof of Theorem~\ref{th:contdecay}}
\label{app:contdecay}

\begin{IEEEproof}
We show by induction on the switching time $t_k$ that the active set $\G_k$ contains less than $q$ active elements and that equations \eqref{eqcont:outputdecay} and
\begin{equation}
 \eltwo{u_{\GD(t)}(t)} \leq e^{-t/\tau} \eltwo{u_{\GD(0)}(0)} + \prn{1-e^{-t/\tau}}\lambda\sqrt{q}, \label{eqcont:stateenergybound}
\end{equation}
hold $\forall t \leq t_k$.

At time $t_0=0$, the theorem hypotheses imply that $\G_0$ contains less than $q$ active elements, and that \eqref{eqcont:outputdecay} and \eqref{eqcont:stateenergybound} hold. 

Next, assume that $\forall t \leq t_k$ the active set contains less than $q$ active elements and that \eqref{eqcont:outputdecay} and \eqref{eqcont:stateenergybound} hold. We start by showing that \eqref{eqcont:outputdecay} holds $\forall t \leq t_{k+1}$. By the induction hypothesis, the active set $\G$ contains less than $q$ active nodes for all $t\leq t_k$, including the current active set $\Gk$ for $t\in[t_k,t_{k+1})$. As a consequence, the inequalities in Lemma \ref{lem:ripcsq} hold with $\G_1 = \G$ and $\G_2 = \Gdag$, for all $t\leq t_{k+1}$.

We compute the following time derivative $\forall t\leq t_{k+1}$:
%
 %
 %
%
\begin{align*}
 \tau \dfe{t}{\dfrac{1}{2}\eltwo{a(t)-\adag(t)}^2}  & = \tau \prn{a(t)-\adag(t)}^T\prn{\dot{a}(t) - \dotadag(t)}\\
 &  = \prn{a(t)-\adag(t)}^T\left(-\PP a(t) + \Phi_{\G}^Ty(t) - \lambda z_{\G} - \tau\dotadag(t) \right) \\
 &  = -\prn{a(t)-\adag(t)}^T\Phi_{\G}^T\Phi_{\prn{\G\cup\Gdag}}\prn{a(t)-\adag(t)} \\
 & \pushright{ + \prn{a(t)-\adag(t)}^T\prn{\Phi_{\G}^T\epsilon(t) - \lambda z_{\G} - \tau\dotadag(t)} }.
\end{align*}
Note that
\begin{align*}
 & -\prn{a(t)-\adag(t)}^T\Phi_{\G}^T\Phi_{\prn{\G\cup\Gdag}}\prn{a(t)-\adag(t)} \\
 & \qquad \qquad = -\eltwo{\Phi_{\G} \prn{a(t)-\adag(t)}}^2 + \prn{a(t)-\adag(t)}^T\Phi_{\G}^T\Phi_{\prn{\G^c\cap\Gdag}}\adag_{\G^c}(t) \\
 & \qquad \qquad \leq -\eltwo{\Phi_{\G} \prn{a(t)-\adag(t)}}^2 + \eltwo{a_{\G}(t)-\adag_{\G}(t)}\eltwo{\Phi_{\G}^T\Phi_{\prn{\G^c\cap\Gdag}}\adag_{\G^c}(t)} \\
 & \qquad \qquad \leq -(1-\delta) \eltwo{a_{\G}(t)-\adag_{\G}(t)}^2  + \delta\eltwo{a_{\G}(t)-\adag_{\G}(t)}\eltwo{\adag_{\G^c}(t)} \\
 & \qquad \qquad = -(1-\delta)\eltwo{a(t)-\adag(t)}^2 +(1-\delta)\eltwo{\adag_{\G^c}(t)}^2 + \delta\eltwo{a_{\G}(t)-\adag_{\G}(t)}\eltwo{\adag_{\G^c}(t)} \\
 & \qquad \qquad \leq -(1-\delta)\eltwo{a(t)-\adag(t)}^2 + \eltwo{\adag_{\G^c}(t)}  \Big( \prn{1-\delta}\eltwo{a(t)-\adag(t)} + \delta\eltwo{a(t)-\adag(t)} \Big)  \\
 & \qquad \qquad = -(1-\delta)\eltwo{a(t)-\adag(t)}^2  + \eltwo{\adag_{\G^c}(t)} \eltwo{a(t)-\adag(t)}.  
 %
\end{align*}
Plugging this back in the expression for the time derivative, we obtain
\begin{align*}
 & \tau \dfe{t}{\dfrac{1}{2}\eltwo{a(t)-\adag(t)}^2} + (1-\delta)\eltwo{a(t)-\adag(t)}^2 \\
 & \qquad \qquad \leq  \eltwo{a(t)-\adag(t)}\eltwo{\adag_{\G^c}(t)}+ \prn{a(t)-\adag(t)}^T\prn{ \Phi_{\G}^T\epsilon(t) - \lambda z_{\G} - \tau\dotadag(t)} \\
 & \qquad \qquad \leq \eltwo{a(t)-\adag(t)}\eltwo{\adag_{\G^c}(t)}  + \eltwo{a(t)-\adag(t)}\eltwo{ \Phi_{\G}^T\epsilon(t) - \lambda z_{\G} - \tau\dotadag(t)} \\
 & \qquad \qquad \leq \eltwo{a(t)-\adag(t)}\eltwo{\adag_{\G^c}(t)} + \eltwo{a(t)-\adag(t)} \Big( \eltwo{\Phi_{\G}^T\epsilon(t)} + \lambda \eltwo{z_{\G}} + \tau\eltwo{\dotadag(t)} \Big) \\
 & \qquad \qquad \leq \eltwo{a(t)-\adag(t)}\prn{\eltwo{\adag_{\G^c}(t)} + \tau\eltwo{\dotadag(t)}} + \eltwo{a(t)-\adag(t)}\prn{ \sqrt{1+\delta}\eltwo{\epsilon(t)} + \lambda \sqrt{q} }\\
 & \qquad  \qquad \leq \eltwo{a(t)-\adag(t)}\prn{\tau\mu + \sigma + \lambda\sqrt{q} }.
\end{align*}
Noting that 
$$\dfrac{d}{dt}\left(\eltwo{x(t)}\right) = \dfrac{\frac{d}{dt}\prn{\eltwo{x(t)}^2}}{2\eltwo{x(t)}}, $$
we obtain
\begin{align*}
 & \dfe{t}{\eltwo{a(t)-\adag(t)}} \leq -(1-\delta)/\tau\eltwo{a(t)-\adag(t)} + 1/\tau\prn{ \sigma + \lambda\sqrt{q} + \tau \mu}.
\end{align*}
Since $\eltwo{a(t)-\adag(t)}$ is continuous, we apply Lemma \ref{lem:odesol} and get $\forall t\leq t_{k+1}$:
\begin{align*}
 \eltwo{a(t)-\adag(t)} & \leq e^{-(1-\delta)t/\tau}\eltwo{a(0)-\adag(0)} + \prn{1-e^{-(1-\delta)t/\tau}} \dfrac{\sigma+\lambda\sqrt{q}+\tau\mu}{1-\delta} .
\end{align*}
This shows that \eqref{eqcont:outputdecay} holds for all $t\leq t_{k+1}$.

We now show that \eqref{eqcont:stateenergybound} holds for all $t\leq t_{k+1}$. Note that, even though the set $\GD(t)$ varies with time, $\eltwo{u_{\GD(t)}(t)}$ is continuous for all $t\geq0$ (as a continuous function of supremum of continuous functions $\abs{u_i(t)}$). Moreover, we can compute the time-derivative $\forall t\leq t_{k+1}$:
\begin{align*}
 \tau\dfrac{d}{dt}\prn{\dfrac{1}{2}\eltwo{u_{\GD}(t)}^2} & = \tau u_{\GD}(t)^T\dot{u}_{\GD}(t) \\
 & = u_{\GD}(t)^T \Big(-u_{\GD}(t) +a_{\GD}(t) - \Phi_{\GD}^T\Phi a(t) + \Phi_{\GD}^Ty(t) \Big)  \\
 %
 & \leq -\eltwo{u_{\GD}(t)}^2 + \eltwo{u_{\GD}(t)}\eltwo{\rho_{\GD}(t)},
\end{align*}
where $\rho_{\GD}(t) = a_{\GD}(t) - \Phi_{\GD}^T\Phi a(t) + \Phi_{\GD}^Ty(t)$. 
We can bound this quantity $\forall t\leq t_{k+1}$:
\begin{align*}
 \eltwo{\rho_{\GD}(t)} & = \eltwo{ a_{\GD}(t) - \Phi_{\GD}^T\Phi a(t) + \Phi_{\GD}^Ty(t)} \\
 & = \eltwo{\adag_{\GD}(t) + \prn{I_{\GD} - \Phi_{\GD}^T\Phi}(a(t) - \adag(t)) + \Phi_{\GD}^T\epsilon(t)} \\
 & \leq \eltwo{\adag(t)} + \eltwo{\Phi_{\GD}^T\epsilon(t)} + \norm{I_{\GD} - \Phi_{\GD}^T\Phi_{\prn{\GD\cup\Gk}}}\eltwo{a(t) - \adag(t)}\\
 & \leq \beta + \sigma + \delta \eltwo{a(t) - \adag(t)},
\end{align*}
where we applied Lemma \ref{lem:ripcsq} with $\G_1 = \GD$ and $\G_2 =\Gdag$, and Lemma \ref{lem:targetenergybound}.
Finally, applying bound \eqref{eqcont:outputdecay} obtained for $\eltwo{a(t) - \adag(t)}$ for all $t\leq t_{k+1}$ and the monotonicity of the exponential, we get:
\begin{align*}
 \eltwo{\rho_{\GD}(t)} & \leq \beta + \sigma + \delta \Big[ e^{-(1-\delta)t/\tau}\eltwo{a(0)-\adag(0)} + \prn{1-e^{-(1-\delta)t/\tau}}D \Big] \\
 & \leq \beta + \sigma + \delta \max\set{\eltwo{a(0)-\adag(0)},D} \\
 & \leq \lambda\sqrt{q},
\end{align*}
where the last inequality comes from the theorem's hypothesis \eqref{eqcont:decaycond}. As a consequence, we obtained that $\forall t\leq t_{k+1}$:
\begin{align*}
\tau\dfrac{d}{dt}\prn{\eltwo{u_{\GD}(t)}} & \leq -\eltwo{u_{\GD}(t)} + \lambda\sqrt{q}.
\end{align*}
Using Lemma \ref{lem:odesol} again yields $\forall t\leq t_{k+1}$
\begin{align*}
\eltwo{u_{\GD}(t)} & \leq e^{-t/\tau}\eltwo{u_{\GD}(0)} + e^{-t/\tau}\dint{0}{t}{ e^{\nu/\tau} \lambda\sqrt{q} d\nu} \\
 & \leq e^{-t/\tau}\eltwo{u_{\GD}(0)} + \prn{1-e^{-t/\tau}}\lambda\sqrt{q},
\end{align*}
which shows that \eqref{eqcont:stateenergybound} holds for all $t\leq t_{k+1}$.

Finally, we prove the last induction hypothesis, that the next active set $\G_{k+1}$ contains less than $q$ indices. Since we proved that \eqref{eqcont:stateenergybound} holds $\forall t\leq t_{k+1}$, together with \eqref{eqcont:stateenergyinit} this implies that:
\begin{align*}
\eltwo{u_{\GD}(t_{k+1})} & \leq e^{-t_{k+1}/\tau}\eltwo{u_{\GD}(0)} +  \prn{1-e^{-t_{k+1}/\tau}}\lambda\sqrt{q} \leq \lambda\sqrt{q}.
\end{align*}
Applying Lemma \ref{lem:inac} shows that the active set $\G_{k+1}$ contains less than $q$ indices and finishes the proof.
\end{IEEEproof}




\bibliographystyle{IEEEbib}
\bibliography{LCAdyn}

\end{document}